\documentclass[11pt, letterpaper,english]{amsart}
\usepackage[left=1in,right=1in,bottom=1in,top=1in]{geometry}
\usepackage{amsmath,amsthm,amssymb}
\usepackage{mathrsfs}
\usepackage[all,cmtip]{xy}
\usepackage{cite}
\usepackage{hyperref}
\usepackage[shortlabels]{enumitem}
\usepackage{graphicx}
\usepackage{lipsum}
\usepackage[usenames,dvipsnames]{color}
\usepackage[normalem]{ulem}

\theoremstyle{plain}

\newtheorem{lem}{Lemma}
\newtheorem{lemma}[lem]{Lemma}
\newtheorem{thm}[lem]{Theorem}
\newtheorem*{thm*}{Theorem}

\newtheorem{cor}[lem]{Corollary}

\newtheorem{que}[lem]{Question}

\theoremstyle{definition}

\numberwithin{equation}{subsection}
\numberwithin{lem}{subsection}

\newcommand{\mathfont}{\mathbb}

\newcommand{\Z}{\mathfont Z}

\DeclareFontFamily{OT1}{rsfs}{}
\DeclareFontShape{OT1}{rsfs}{n}{it}{<-> rsfs10}{}
\DeclareMathAlphabet{\mathscr}{OT1}{rsfs}{n}{it}

\DeclareMathOperator{\Ker}{ker}

\DeclareMathOperator{\et}{\acute{e}tale}

\DeclareMathOperator{\Spec}{Spec}

\usepackage{tikz-cd}

\title[Local-to-Global Extensions]{Local-to-Global Extensions to Wildly Ramified Covers of Curves}
\date{\today}

\author{Renee Bell}
\email{rhbell@math.upenn.edu}

\begin{document} 

\maketitle

\begin{abstract}
Given a Galois cover of curves $X \to Y$ with Galois group $G$ which ramified only at $x$, restriction to the punctured formal neighborhood of $x$ induces a Galois extension of Laurent series rings $k((u))/k((t))$. If we fix a base curve $Y$, we can ask when a Galois extension of Laurent series rings comes from a global cover of $Y$ in this way. Harbater proved that over a separably closed field, every Laurent series extension comes from a global cover for any base curve if $G$ is a $p$-group, and gave a condition for the uniqueness of such an extension. Using a generalization of Artin--Schreier theory to non-abelian $p$-groups which gives an explicit description of $G$-covers, we fully characterize the curves $Y$ for which this extension property holds and for which it is unique up to isomorphism, but over a more general ground field. 
\end{abstract}

\section{Introduction}

Throughout this paper, $k$ is an arbitrary field of characteristic $p$, and $G$ is a finite $p$-group.

Let $Y$ be a smooth proper curve over $k$ and $y \in Y (k)$. We define a ``$y$-ramified $G$-cover of $Y$" to
be a Galois cover of curves $q : X \to Y$ with Galois group $G$, 
which is unramified
on $Y' := Y - \{y\}$. By the Cohen structure theorem, we can choose uniformizers $t$, $u_i$ such that 
$\prod_{x_i \in q^{-1}(y)}\widehat{\mathcal{O}_{X,x_i}}
\cong
\prod_{x_i \in q^{-1}(y)} k [[u_i]]
$
 and $\widehat{\mathcal{O}_{Y,y}} \cong k[[t]]$.
 After localization, we obtain a $G$-Galois \'{e}tale algebra $L:= \prod_{x_i \in q^{-1}(y)} k ((u_i))$ over $k((t))$. We say that $L$ arises from the $G$-action on $X$. 
 
 Thus, for each curve $Y$ and point $y \in Y$ and $p$-group $G$, we obtain a functor
\[
\psi_{Y,y,G}: 
 \Bigg\{
\begin{gathered}
  \textrm{$y$-ramified}\\
  G\textrm{ covers}\\
  \textrm{ of } Y
\end{gathered}
  \Bigg\}
  \to
   \Bigg\{
   \begin{gathered}
  \textrm{Galois \'{e}tale algebras }\\
  \textrm{over }k((t))\textrm{ with Galois}\\
  \textrm{group a } G
\end{gathered}
  \Bigg\}
\]

Understanding the functor $\psi_{Y,y,G}$ allows us to use the geometry of Galois covers of curves to classify automorphisms of $k[[t]]$ as in \cite{BCPS}.
Conversely, it allows us to use extensions of $k((t))$ in order to classify filtrations of ramification groups of Galois covers of curves with Galois group a $p$-group as noted in the survey \cite{HOPK}.

Questions about $\psi_{Y,y,G}$ can be approached by turning to \'{e}tale cohomology. Throughout this
paper, for a scheme $S$ and a (not necessarily abelian) group $G$, we denote by $H^1(S,G)$ the \v{C}ech cohomology $\check{H}^1_{et} (S,G)$ of the constant sheaf of groups with coefficients in $G$ with respect to the \'{e}tale site on $X$; this cohomology set parameterizes principal $G$-bundles on $X$ \cite{milneLEC}.
The map from
the ring of regular functions on $Y'$ into the Laurent series field $\mathcal{O}(Y') \to k((t))$ coming from the Cohen Structure Theorem induces a map
 $\Spec k((t)) \to Y'$ which we can think of as inclusion of the formal deleted neighborhood around $y$ 
into $Y'$. Hence, we obtain a map $H^1(Y', G) \to H^1(k((t)), G)$ which we denote by $\Psi_{Y,y,G}$. We note that $\Psi_{Y,y,G}$ is induced from $\psi_{Y,y,G}$ by passing to isomorphism classes.

We pose some basic questions about $\Psi_{Y,y,G}$.

\begin{que}\label{surjq}
When is  $\Psi_{Y,y,G}$ surjective?
\end{que}

This is equivalent to asking when every $G$-Galois extension of $k((t))$ extends to a global, $y$-ramified Galois cover of $Y$. 
In \cite{harbater},
Harbater showed that if the ground field $k$ is algebraically closed, then $\Psi_{Y,y,G}$ is surjective for any $p$-group $G$. In this paper, we provide an answer to Question~\ref{surjq}
over a more general field $k$, not necessarily algebraically or even separably closed, in the following theorem.

Notation: for any ring $R$ of characteristic $p$, let $\wp: R\to R$ denote the Artin--Schreier map $f \mapsto f^p - f$.

\begin{thm}\label{surj}
Let $G$ be a nontrivial finite $p$-group. Then the following are equivalent:
\begin{enumerate}
    \item The equality $k((t)) = \wp(k((t))) + \mathcal{O}(Y')$ holds.
    \item The map $\Psi_{Y,y,G}$ is surjective.
    \item The map $\Psi_{Y,y,\Z/p\Z}$ is surjective.
\end{enumerate}
\end{thm}

We can also ask when any lift of an extension of $k((t))$ to a global Galois cover of $Y$ is unique up to isomorphism.

\begin{que}\label{injq}
When is  $\Psi_{Y,y,G}$ injective?
\end{que}

An answer to this over $k$ algebraically closed was given as well by Harbater in \cite{harbater}.
In fact, he calculates the size of the fiber of $\Psi_{Y,y,G}$ as $p^r$, where $r$ is the $p$-rank of $Y$. We extend the answer to Question~\ref{injq}
to a more general field $k$, which may not be algebraically closed.

\begin{thm}\label{inj}
Let $G$ be a nontrivial finite $p$-group. Then the following are equivalent:
\begin{enumerate}
    \item The equality $\wp(k((t))) \cap \mathcal{O}(Y') = \wp(\mathcal{O}(Y'))$ holds.
    \item The map $\Psi_{Y,y,G}$ is injective.
    \item The map $\Psi_{Y,y,\Z/p\Z}$ is injective.
\end{enumerate}
\end{thm}

Combining our answers to Questions~\ref{surjq} and \ref{injq}, without assuming that the base field is algebraically or even separably closed, gives a criterion on $Y$ for $\psi_{Y,y,G}$ to be an equivalence of categories. This generalizes the result of Katz in \cite{katz}
which states that over any field of characteristic $p$, the functor $\psi_{\mathbb{P}^1,\infty}$ is an equivalence of categories. Curves satisfying the criteria of Theorems~\ref{surj} and 
\ref{inj}
are particularly useful for relating the geometry of the curve and its covers to properties of $k((t))$ and its extensions. 
In Section 5
of this paper, we give another explicit example of a class of such curves.

Our proofs of these theorems use new and more explicit methods. Proofs in previous work, as in \cite{katz},
have reduced the problem to the case in which $G$ is abelian. In this case, one can use the vanishing of certain $H^2$ groups or a characterization of abelian $p$-group field extensions using Witt vector theory, as noted in \cite{HOPK}.
In this paper, we describe and work with an explicit characterization of $G$-Galois \'{e}tale algebras for $G$ not necessarily abelian. This characterization, which we will call the Inaba classification, is a generalization of a theorem of Inaba in \cite{inaba}
which extends Artin--Schreier--Witt theory to nonabelian Galois \'{e}tale algebras.

We introduce some notation: let $U_n(R)$ denote the group of upper triangular $n \times n$ matrices with entries in $R$ such that all diagonal entries are $1$, and.
%raise to the p later

In order to be more explicit about the structure of these extensions, we put an order on the indices of entries of matrices in $U_n(R)$. The indices $(i,j): 1\le i<j\le n$ are ordered by 
    lexicographical order on $(j-i,j)$; that is, they are ordered by going down along consecutive diagonals.

\begin{thm}\label{inaba}
Let G be a finite p-group, and fix an injective homomorphism $\Lambda : G \to U_n(\mathbb{F}_p)$ for some suitable $n$. Let $R$ be a ring of characteristic $p$ such that $\Spec R$ is connected, and let $L/R$ be a Galois  \'{e}tale algebra with Galois group $G$.
\begin{enumerate}
    \item The $R$-algebra $L$ is generated by elements $a_{ij} \in L$ for $1 \leq i < j \leq n$ such that the unipotent matrix $A := (a_{ij})$ satisfies $A^{(p)} = MA$ for some $M \in U_n(R)$. We also have that for $\sigma \in G$, ${}^{\sigma}A = A\Lambda(\sigma)$, where $\sigma$ acts entry-wise on A. We denote this situation by $L = R[A]$ and say that $L$ is of type $M$.
    \item Given two algebras $L,L' \in H^1(\Spec R,G)$, if we choose $(A,M)$ for $L$ and $(A',M')$ for $L'$, then $L,L'$ are isomorphic if and only if $M = C^{(p)}M'C^{-1}$ for some $C \in U_n(R)$.
    \item In addition, $L/R$ can be decomposed as a composition of subextensions 
\[
    L = L_{1,n}/L_{2,n}/L_{1,n-1}/ ... /L_{2,3}/ L_{1,2}/R \label{extension} \tag{*}
\]
    where 
    $L_{i,j} = R[\{ a_{\ell m} : (\ell,m) \leq (i,j) \}]$. 
    And let $L_{<i,j}$ denote $R[\{ a_{\ell m} : (\ell,m) < (i,j) \}]$. 
    Then for each $(i,j)$, either $a_{ij} \in L_{<i,j}$ or $L_{i,j}\cong L_{<i,j}[x]/(x^p-x-(\sum_{i<k<j} m_{ik}a_{kj})- m_{ij})$. The latter is the case if and only if there is an element $g\in G$ such that $\Lambda(g)_{i,j}\ne 0$ but $\Lambda(g)_{l,m}=0$ whenever $(l,m)<(i,j)$. 
\end{enumerate}
\end{thm}

Note that the above is exactly the Artin--Schreier characterization of $\Z/p\Z$-Galois \'{e}tale algebras when $G = \Z/p\Z \cong U_2(\mathbb{F}_p)$.

Our proof of the Inaba classification uses modern methods of \'{e}tale cohomology and hence provides an alternative proof of Inaba's result while also generalizing the result. Toward this end, we present several useful results in nonabelian \'{e}tale cohomology in Section 2.

\section{Exact Sequences for Nonabelian Cech Cohomology}

In order to prove our results in the setting where $G$ is nonabelian, we present some results about exact sequences of \v{C}ech cohomology.

If $\mathcal{G}$ is a sheaf of abelian groups, we define $\check{H}^i(X,\mathcal{G})$ 
to be the $i$th \v{C}ech cohomology group with respect to the \'{e}tale site, 
and if $X$ is affine, then $\check{H}^i(X,\mathcal{G})$ is isomorphic to $H^i_{et}(X,\mathcal{G})$ by Theorem 10.2 of \cite{milneLEC}.

But if $\mathcal{G}$ is a sheaf of groups which are nonabelian, only $\check{H}^0(X,\mathcal{G})$ and $\check{H}^1(X,\mathcal{G})$ are defined, and the latter may not be a group.

\subsection{Nonabelian \v{C}ech cohomology}

First, we recall definitions of these cohomology sets and describe natural maps between them, following \cite{milneLEC}.

Let $\mathcal{G}$ be a sheaf of groups on $X_{et}$
and 
let $\mathcal{U}= (U_i \to X)_{i \in I}$ be an \'{e}tale covering of $X$. 
Let $U_{{i_1} \cdots {i_k}}$ denote $U_{{i_1}}\times_X  \cdots \times_X U_{{i_k}}$.
A 1-cocycle for $U$ with values in $\mathcal{G}$ is a family $(g_{ij})_{(i,j) \in I \times I}$ with $g_{ij} \in \mathcal{G}(U_{ij})$ such that 
\[
(g_{ij}|U_{ijk}) \cdot (g_{jk}|U_{ijk}) = g_{ik}|U_{ijk} \text{ for all } i,j,k.
\]
We define an equivalence relation $\sim $ as follows; for two cocycles $g = (g_{ij})$ and $g'=(g'_{ij})$, we write $g \sim g'$ if there is a family $(\gamma_i)_{i \in I}$ with $\gamma_i \in \mathcal{G}(U_i)$ such that 
\[
g'_{ij} = (\gamma_i | U_{ij}) \cdot g_{ij} \cdot (\gamma_j | U_{ij})^{-1}
\]
for all $i,j$. Then $\check{H}^1(X_{et},\mathcal{G})$ is defined to be the limit over all \'{e}tale coverings of $X$ of 1-cocycles modulo $\sim$.

As mentioned before and in \cite{milneLEC}, $\check{H}^1(X_{et},\mathcal{G})$ parameterizes principal $\mathcal{G}$-bundles (when $\mathcal{G}$ is a constant sheaf with group $G$, these are principal $G$-bundles), but it is not necessarily a group. However, $\check{H}^1(X_{et},\mathcal{G})$ has a distinguished element $g_{ij} = 1$ for all $i,j$. 
If $\phi$ is a map to $\check{H}^1(X_{et},\mathcal{G})$, we define 
$\Ker \phi$ to be the preimage of this distinguished element.

\subsection{Partial delta functors}

In the setting of \v{C}ech cohomology, we also have partial delta functors connecting cohomology sets, which we define below.

If 
\[
1 \to \mathcal{G}' \to \mathcal{G} \to \mathcal{G}'' \to 0
\]
is an exact sequence of sheaves of groups, we can define a map 
\[
\delta_0: \mathcal{G}''(X) \to \check{H}^1(X_{et},\mathcal{G}')
\]
as follows.
Let  $c$ be an element of  $\mathcal{G}''(X)$. Then since $\pi: \mathcal{G}\to \mathcal{G}''$ is locally surjective, there exists an \'{e}tale cover $U= (U_i)_{i \in I}$ of $X$, and elements $b_i \in \mathcal{G}(U_i)$ such that $ \pi(b_i) = c|U_i$. We define 
\[
\delta_0(c) = (b_i^{-1}b_j)_{ij}.
\]
One can check that $\delta_0(c)$ is in $\check{H}^1(X_{et},\mathcal{G}')$ and that it is well-defined. 
We also note that $\delta_0$ can be defined in terms of the characterization of $H^1(X,\mathcal{G}')$ as principal $\mathcal{G}'$-bundles by mapping $c\in \mathcal{G}''(X)$ to the subsheaf of $\mathcal{G}$ defined on each open as the preimage of $c$.

When $\mathcal{G}'$ is a sheaf of abelian groups such that $\mathcal{G}'(U)$ is central in $\mathcal{G}(U)$ for all $U$ in an \'{e}tale cover, we can define a map 
\[
\delta_1: 
\check{H}^1(X_{et},\mathcal{G}'') \to \check{H}^2(X_{et},\mathcal{G}')
\]
as follows. 

Let $\mathcal{U} = (U_i \to X)_{i \in I}$ be an \'{e}tale covering of $X$, and let $(c_{ij})_{i \in I}$ be an element of $\check{H}^1(X_{et},\mathcal{G}'')$,
with $c_{ij} \in \mathcal{G}''(U_{ij})$.
Let $b_{ij} \in \mathcal{G}(U_{ij})$ such that $\pi(b_{ij}) = c_{ij}$.
%LATER is this a problem with LOCAL surjectivity? help
Then, from the cocycle condition on $(c_{ij})$, we have that for each triple $i, j,k$, there exists $h_{ijk} \in \mathcal{G}'(U_{ijk})$ such that 
$
h_{ijk} = b_{ik}^{-1}b_{ij}b_{jk}.
$
We define
\[
\delta_1((c_{ij})_{i,j}) = (h_{ijk})_{i,j,k}.
\]
The map $\delta_1$ takes values in $\check{H}^2(X_{et},\mathcal{G}')$ and is well-defined, as can be checked directly and is proved in Section 4.2 of Chapitre III of \cite{giraud}.

\subsection{Long exact sequences}

We now have maps that we can arrange into a sequence which we will show to be ``exact," which we define in the usual sense but with kernels of \textit{pointed sets}. 
Additionally, we may have group actions on cohomology sets which allow for a classification of the fibers; we describe these group actions below.

Let $b \in \mathcal{G}(X)$ and $c \in (\mathcal{G}/\mathcal{H})(X)$. So there is an \'{e}tale cover $\{U_i\}_{i \in I}$ and elements $b_i \in \mathcal{G}(U_i)$ such that $c|U_i$ is the image of $b_i$ for all $i$. Then $b\cdot c\in (\mathcal{G}/\mathcal{H})(X)$ denote the section whose restriction to $U_i$ is the image of $b \cdot b_i$ for each $i$.

If $\mathcal{H}(U)$ is a normal subgroup of $\mathcal{G}(U)$ for all \'{e}tale opens $U$, then we have an action of $(\mathcal{G}/\mathcal{H})(X)$ on $H^1(X, \mathcal{H})$ as follows.
If $c \in (\mathcal{G}/\mathcal{H})(X)$, there is a family $(b_i)_i$ with $b_i \in \mathcal{G}(U_i)$ such that the image of $b_i$ in $\mathcal{G}(U_i)/\mathcal{H}(U_i)$ is $c|U_i$. 
If $a = (a_{ij})$ is an element of $H^1(X,\mathcal{H})$, 
then $c \cdot a := (b_i a_{ij} b_j^{-1})_{ij}$.

If $\mathcal{H}(U)$ is a central subgroup of $\mathcal{G}(U)$ for all \'{e}tale open sets $U$, then we have an action of the abelian group $H^1(X,\mathcal{H})$ on $H^1(X,\mathcal{G})$ as follows. 
If $a = (a_{ij})$ is an element of $H^1(X,\mathcal{H})$ and $b = (b_{ij})$ is an element of $H^1(X,\mathcal{G})$,
then $a\cdot b :=(a_{ij}b_{ij})$. We summarize some observations about fibers of cohomology maps with this in mind.

\begin{lem}\label{longexact}
Let $X$ be a connected scheme, and let $\mathcal{G}$ be a sheaf of (not necessarily abelian) groups, and $\mathcal{H}$ a subsheaf of groups.
\begin{enumerate}
    \item The \v{C}ech cohomology sequence 
\[
1 \to \mathcal{H}(X) \to 
\mathcal{G}(X) \to
(\mathcal{G}/\mathcal{H})(X) \xrightarrow{\delta_0} 
H^1(X, \mathcal{H}) \xrightarrow{\iota}
H^1(X, \mathcal{G}) \xrightarrow{\pi}
H^1(X, \mathcal{G}/\mathcal{H})
\]
is exact, and for elements $c,c' \in (\mathcal{G}/\mathcal{H})(X)$, $\delta_0 (c) = \delta_0(\tilde{c})$ if and only if there exists $\beta \in \mathcal{G}(X)$ such that $c = \beta \tilde{c}$. 

\item If, in addition, $\mathcal{H}(U)$ is a normal subgroup of $\mathcal{G}(U)$ for all \'{e}tale open sets $U$, 
 two elements of $H^1(X, \mathcal{H})$ have the same image in $H^1(X, \mathcal{G})$ if and only if they are in the same $(\mathcal{G}/\mathcal{H})(X)$-orbit.

\item If, in addition, $\mathcal{H}(U)$ is a central subgroup of $\mathcal{G}(U)$ for all \'{e}tale open sets $U$,
then the extended sequence
\[
1 \to \mathcal{H}(X) \to 
\mathcal{G}(X) \to
(\mathcal{G}/\mathcal{H})(X) \to 
H^1(X, \mathcal{H}) \xrightarrow{\iota}
H^1(X, \mathcal{G}) \xrightarrow{\pi}
H^1(X, \mathcal{G}/\mathcal{H}) \xrightarrow{\delta_1} 
H^2(X,\mathcal{H})
\]
is exact, and 
then two elements of $H^1(X, \mathcal{G})$ have the same image in $H^1(X,\mathcal{G}/\mathcal{H})$ if and only if they are in the same $H^1(X,\mathcal{H})$-orbit. 
\end{enumerate}
\end{lem}

\begin{proof}
We prove part (3); its proof is similar to the proofs of parts (1) and (2), which follow from Propositions~3.3.3 and 3.4.5 of Chapitre III of \cite{giraud}.

Proof of part (3): an element $(c_{ij})_{i,j}$ is in the kernel of the map $H^1(X, \mathcal{G}/\mathcal{H}) \to H^2(X, \mathcal{H})$ if and only if representatives $b_{ij}$ for $c_{ij}$ in $\mathcal{G}(U_{ij})$
satisfy the cocycle condition, which implies exactness.

Now we consider the fibers of $\pi$, where $\pi$ denotes the map $H^1(X, \mathcal{G}) \to H^1(X, \mathcal{G}/\mathcal{H})$. 
Let $b, b'$ be elements of $H^1(X,\mathcal{G})$ represented by cocycles $(b_{ij})_{i,j}$ and $(b'_{ij})_{i,j}$ respectively, and let $c_{ij}$ denote the image of $b_{ij}$ in $(\mathcal{G}/\mathcal{H})(U_{ij})$ (and similarly for $c'_{ij}$).
Suppose that  $\pi(b) = \pi(b')$, so there exist $\gamma_i \in (\mathcal{G}/\mathcal{H})(U_i)$ such that $c_{ij} = \gamma_i c'_{ij} \gamma_j^{-1}$ for all $i,j$. 
Now let $\beta_i \in \mathcal{G}(U_i)$ such that the image of $\beta_i$ in $\mathcal{G}/\mathcal{H}(U_i)$ is $\gamma_i$.
Then $a_{ij} := b_{ij}^{-1} \beta_i b'_{ij} \beta_j^{-1}$ is an element of $\mathcal{H}(U_{ij})$.
Since $\mathcal{H}(U_{ij})$ is central in $\mathcal{G}(U_{ij})$, we also have that $a_{ij} = \beta_i b_{ij}' \beta_{j}^{-1} b_{ij}^{-1}$. Centrality of $\mathcal{H}$ together with the cocycle condition on $(b_{ij})_{i,j}$ and $(b'_{ij})_{i,j}$ give that $(a_{ij})_{i,j}$ is a cocycle, and so $b$ and $b'$ are in the same $H^1(X, \mathcal{H})$-orbit. The converse is straightforward.

\end{proof}

We close this section with a useful lemma about maps of cohomology induced by subgroups of $p$-groups.

\begin{lem}\label{subgroup}
Let $R$ be a ring of characteristic $p$ such that $\Spec R$ is connected.
If $G'$ is a $p$-group and $G$ is a subgroup of $G'$, then the induced map $\rho: H^1(\Spec R,G) \to H^1(\Spec R,G')$ is injective.
\end{lem}

\begin{proof}
We first prove the statement for $G$ central in $G'$.
%with $G \cong \Z/p\Z$.
In this case, the action of $G'/G$ on $H^1(\Spec R, G)$ is trivial, so by Lemma~\ref{longexact}, $\iota$ is injective.

We now proceed by induction on $|G'|$, doing a diagram chase of pointed sets.
Let $G'_1$ be a nontrivial central subgroup of $G'$ which is isomorphic to $\Z/p\Z$ (every $p$-group has a non-trivial center, as guaranteed by the class equation),
and let $G_1 = G'_1 \cap G$. By Lemma~\ref{longexact} and Lemma 1.4.3 of \cite{katz}
(the latter stating that $H^2(\Spec R, \Z/p\Z)= 0$),
the rows of the following commutative diagram are exact.
\[
\xymatrix{
1 \ar[r]& 
H^1(\Spec R,G_1)\ar[r]^{\iota} \ar[d]^{\rho_{1}} &
H^1(\Spec R,G) \ar[r]^{\pi} \ar[d]^{\rho} & 
H^1(\Spec R,G/G_1) \ar[r] \ar[d]^{\rho_2} & 
1\\
1 \ar[r] & 
H^1(\Spec R,G'_1) \ar[r]^{\iota'}  & 
H^1(\Spec R,G') \ar[r]^{\pi'}  & 
H^1(\Spec R,G'/G'_1) \ar[r]  & 
1
}
\]
By induction, $\rho_1$ and $\rho_2$ are injective, and by the case in the beginning of this proof, $\iota$ and $\iota'$ are injective.
Now suppose that $b_1,b_2 \in H^1(\Spec R, G)$ and $\rho(b_1) = \rho(b_2)$.
Then $\rho_2\circ \pi(b_1)=\rho_2\circ \pi(b_2)$, so $\pi(b_1)=\pi(b_2)$.
By Lemma~\ref{longexact}, there exists $a \in H^1(\Spec R, G_1)$ such that $a \cdot b_1 = b_2$.
Then $\rho_1(a)\rho(b_1) = \rho(b_2)=\rho(b_1)$, so $\rho_1(a) = 1$ by injectivity of $\iota'$, and so $a = 1$ by injectivity of $\rho_1$. Thus $b_1 = b_2$ and $\rho$ is injective.
\end{proof}

\section{Inaba Classification of p-group Covers}
We now provide generalizations of Artin--Schreier theory to non-abelian groups. 
We first establish some notation. 
For an $\mathbb{F}_p$-algebra $R$, we define $U_n(R)$ to be the group of upper triangular $n\times n$ matrices with coefficients in $R$, such that all the diagonal entries are 1, and the group action is matrix multiplication. We also denote by $X$ the upper triangular matrix of indeterminates, where the $(i,j)$ entry is the indeterminate $x_{ij}$ for $j-i>0$ and the diagonal entries are 1. 
For a characteristic $p$ ring $R$ and a matrix $C \in U_n(R)$, we denote by $C^{(p)}$ the matrix obtained by raising each entry of $C$ to the $p$th power (different from matrix multiplication of $C$ with itself
$p$ times).
For a matrix $M \in U_n(R)$, we denote by $L_M$ the $R$-algebra $R[X]/(X^{(p)} - MX)$, by which we mean
the $R$-algebra generated by the indeterminate entries of $X$, modulo the relations
coming from the matrix equation $X^{(p)} = MX$. This has a $U_n(\mathbb{F}_p)$-action given by
$X \mapsto X \cdot g$ for $g \in U_n(\mathbb{F}_p)$ (with indeterminates mapping to the corresponding entry
of the matrix $X \cdot g$). 
Lastly, we say that two matrices $M,M' \in U_n(R)$ are $p$-equivalent over $R$ if there exists $C \in U_n(R)$ such that $M = C^{(p)}M' C^{-1}$.

\begin{lemma}\label{mult}
Let $R$ be a ring of characteristic $p$ such that Spec $R$ is connected. 
Then the finite Galois \'{e}tale algebras over Spec $R$ with Galois group $U_n(\mathbb{F}_p)$ are the algebras $R[X]/(X^{(p)} = MX)$ where $M$ ranges over all matrices in $U_n(R)$, and the Galois action is given by matrix multiplication $X \mapsto X\cdot g$. 
Two such Galois algebras defined by matrices $M, M'$ are isomorphic as $R$-algebras with $U_n(\mathbb{F}_p)$-action if and only if $M = C^{(p)}M'C^{-1}$ for some $C \in U_n(R)$.
\end{lemma}

\begin{proof}
Let $U_n$ be the $\mathbb{F}_p$-group scheme representing the functor which sends a ring $A$ to $U_n(A)$, and let $U_n(\mathbb{F}_p)$ be the constant group scheme. We have an exact sequence
\[
1 \to U_n(\mathbb{F}_p) \to U_n \xrightarrow{\mathcal{L}} U_n \to 1
\]
where $\mathcal{L}$ is the morphism (which is not a group homomorphism) $B \mapsto B^{(p)} B^{-1}$.
By Lang's theorem \cite{langsurj}
$\mathcal{L}$ is surjective and identifies $U_n/(U_n(\mathbb{F}_p))$ with $U_n$. 
Since Spec $R$ is connected, we know that $H^0(\textrm{Spec }R, U_n(\mathbb{F}_p)) = U_n(\mathbb{F}_p)$, so 
by Lemma \ref{longexact}
we have an exact sequence of pointed sets 
\[
1 \to U_n(\mathbb{F}_p) \to 
U_n(R) \xrightarrow{\mathcal{L}} 
U_n(R) \xrightarrow{\delta}
H^1(\textrm{Spec } R, U_n(\mathbb{F}_p)) \to 
H^1(\textrm{Spec } R, U_n)
\]
where $\delta$ sends a matrix $M \in U_n(R)$ to the principal $U_n(\mathbb{F}_p)$-bundle given by $\mathcal{L}^{-1}(M)$ on each \'{e}tale open of $\Spec R$. That is, $\delta(M)$ is exactly the \'{e}tale algebra $R[X]/(X^{(p)}=MX)$.

Since $H^1(X,\mathcal{O}_X) = 1$ for $X$ affine and $U_n$ has a composition series whose factors
are $\mathbb{G}_a$, 
we see by induction that $H^1(\Spec R,U_n)=1$,
so the map $U_n(R) \to H^1(\Spec R, U_n(\mathbb{F}_p))$ is surjective, and every element of $H^1(\Spec R, U_n(\mathbb{F}_p))$ can be represented by a $U_n(\mathbb{F}_p)$-algebra of the form $L_M$.
By Lemma \ref{longexact}, $L_M \cong L_{M'}$ if and only if there exists a matrix $C \in U_n(R)$ such that $M' = C * M$, where $*$ is the left action of $U_n(R)$ on $U_n/U_n(R)$.
Let $N$ be an element of $U_n(S)$ for $S$ an \'{e}tale $R$-algebra such that $\mathcal{L}(N) = M$.
Then $C * M = \mathcal{L}(CN) = C^{(p)} M C^{-1}$. This yields the result.
\end{proof}

Now we look at $G$ a general $p$-group and fix an embedding $\Lambda : G \to U_n(\mathbb{F}_p)$ (such an embedding is guaranteed by Proposition~2.4.12 of \cite{springer}.

\subsection{Proof of Theorem \ref{inaba}}
Part (1):
First, we note that the inclusion $\Lambda: G \to U_n(\mathbb{F}_p)$ induces a map 
\[
H^1(\Spec R, G) \to H^1(\Spec R,U_n(\mathbb{F}_p))
\]
sending
\[
L \mapsto \prod_{G \backslash U_n(\mathbb{F}_p)} L =: \tilde{L}
\]
with the following left $U_n(\mathbb{F}_p)$-action.
Let $u_1,...,u_r$ be coset representatives for $G \backslash U_n(\mathbb{F}_p)$, with $u_1 = e$ being the identity element. Then we can write any element of $\prod_{G \backslash U_n(\mathbb{F}_p)} L$ as $(\ell_i)_{i=1}^r$ with $\ell_i \in L$. For each $u \in U_n(\mathbb{F}_p)$, there exist $g_i \in G$ such that $u_i u = g_i u_{j(i)}$, where $j(i)$ is the index of the coset $u_i u$. 
Then $u\cdot (\ell_i)_{i=1}^r = (^{g_{j^{-1}(i)}} \ell_{j^{-1}(i)})_{i=1}^r$.

Now consider the map $\pi: \tilde{L} \to L$ which is projection onto the first component, and note that 
for $g \in G$, the first coordinate of $g\cdot (\ell_i)_{i=1}^r$ is $^{g} \ell_1$,
so $\pi$ is a map of \'{e}tale $G$-algebras. 
But by Lemma~\ref{mult}, $\tilde{L} \cong R[X] / (X^{(p)} = MX)$ as $G$-algebras, so the surjection $\pi$ expresses $L$ as $R[A]$, where $A$ is the matrix with $ij$-coordinate $\pi(x_{ij})$. And since $\pi$ is compatible with the action of $G$, the original $G$-action on $L$ agrees with the action coming from matrix multiplication by $\Lambda(G)$.

Part (2): following \cite{inaba}, we show that the matrix $M$ is unique up to $p$-equivalence.

Suppose that $L = R[\{ a_{ij} \}]$ and  $A:= (a_{ij})$ satisfies $A^{(p)} = MA$
for some $M \in U_n(R)$ such that the action of $G$ is given by ${}^\sigma A = A \cdot \sigma$. Suppose that $L$ can similarly be generated by the entries of $B$ with $B^{(p)} = \tilde{M}B$
for some $\tilde{M} \in U_n(R)$. 
Then $B A^{-1} \in U_n(R)$ since $B A^{-1}$ is fixed by $G$,
so $B = CA$ for some $C \in U_n(R)$,
and thus $\tilde{M} = C^{(p)}MC^{-1}$.

Conversely, suppose $L,L'$ are $G$-Galois $R$-algebras of type $M$ and $M'$, respectively, and that $M$ is $p$-equivalent to $M'$. 
By the above arguments, $L$ maps to $L_M$ in $H^1(\Spec R,U_n(\mathbb{F}_p))$ and $L'$ maps to $L_{M'}$, and by Lemma~\ref{mult} $L_M$ and $L_{M'}$ are isomorphic. Since the map $H^1(\Spec R, G) \to H^1(\Spec R, U_n(\mathbb{F}_p))$ is injective by Lemma~\ref{subgroup}, $L$ must be isomorphic to $L'$. 

Part (3): First note that, since $L/R$ has degree a power of $p$, each of the subextensions in $(*)$ must have degree a power of $p$. Let $1\le i<j\le n$. The equation $A^{(p)}=MA$ at the $(i,j)$ coordinate gives $a_{ij}^p = a_{ij} + (\sum_{i<k<j} m_{ik}a_{kj}) + m_{ij}$. Since $L_{i,j}=L_{<i,j}[a_{ij}]$, the extension $L_{i,j}/L_{<i,j}$ must have degree $1$ or $p$, giving the first claim.

For the second claim, note the Galois group of $L$ over $L_{i,j}$ is $G \cap N_{i,j}$, where $N_{i,j}$ is the subgroup of matrices in $U_n(\mathbb{F}_p)$ with $(\ell,m)$-entry equal to zero for $(\ell,m)\le (i,j)$ (for entries above the main diagonal). 
We similarly have that $G\cap N_{<i,j}$ (matrices with $(\ell,m)$-entry equal to zero for $(\ell,m)< (i,j))$ is the Galois group of $L$ over $L_{<i,j}$. Hence $(G\cap N_{<i,j})/(G\cap N_{i,j})$ is the Galois group of $L_{i,j}/L_{<i,j}$, giving the second claim.
\qed

\section{Main Proofs}

We now show that properties of the map $\Psi_{Y,y,G}$ can be checked for $G = \Z/p\Z$, or equivalently that it suffices to know the behavior of $\wp$ on $k((t))$ and $\mathcal{O}(Y')$. We begin with a lemma about the structure of $U_n(R)$.

\begin{lem}\label{triangle}
Let $R$ be a ring of characteristic $p$. Suppose $M = (m_{ij} )$, $M' = (m'_{ij} ) \in U_n (R)$ are $p$-equivalent, so $M = B^{(p)}M'B^{-1}$ for some $B = (b_{ij}) \in U_n(R)$.
Then for each pair $i,j,$ there exists an element $c$ of the $\Z$-algebra generated by $\{ m_{i' j'} , b_{i' j'} , m'_{i' j '} | i ' - j ' < i - j \}$ such that $m_{ij} = \wp(b_{ij} ) + m'_{ij} + c$ . 
That is, $m_{ij} = \wp(b_{ij} ) + m'_{ij}$ modulo the elements on the lower diagonals.
\end{lem}

\begin{proof}
Consider two matrices $W = (w_{ij})$, $Z = (z_{ij}) \in U_n(R)$. We compute that the $(i,j)$ entry of $WZ$ is 
\[
\sum_{k=1}^n w_{ik}z_{kj}
=
z_{ij} + w_{ij} + \sum_{i<k<j} w_{ik}z_{kj}.
\]
Applying this to both sides of the equation $MB = B^{(p)}M'$ yields the result.
\end{proof}

\subsection{Proof of Theorem~\ref{inj}}
We first show that (1) implies (2), so suppose $ \wp(k((t)))\cap \mathcal{O}(Y') = \wp (\mathcal{O}(Y'))$.
Let $\Spec L, \Spec L'$ be two \'{e}tale $G$-covers of $Y'$.
By Theorem~\ref{inaba}
$L=R[A]$ and $L' = R[A']$ with $A^{(p)} = MA$ and $A'^{(p)} = M'A'$ for some $M,M' \in U_n(\mathbb{F}_p)$, and since they are isomorphic over $k((t))$, we have that $M = B^{(p)}M' B^{-1}$ for some $B \in U_n(k((t)))$.
Suppose for contradiction that there exists an entry $b_{ij}$
of $B$ not in $\mathcal{O}(Y')$, chosen such that all entries on lower diagonals are in $\mathcal{O}(Y')$. Then by Lemma~\ref{triangle}
, $\wp(b_{ij}) = m_{ij} - m'_{ij} + c$ where $c$ is a polynomial in the entries of lower diagonals of $B, M$, and $M'$. 
Then $\wp(b_{ij} ) \in \mathcal{O}(Y')$, and by (1), $b_{ij} \in \mathcal{O}(Y')$, a contradiction.

Next we show that (2) implies (3), so assume that $\Psi_{Y,y,G}$ is injective. Since $G$ is a nontrivial $p$-group, it has a nontrivial subgroup $H$ which is isomorphic to $\Z/p\Z$,
so by Lemma~\ref{subgroup}, we have a commutative square 

\[
\xymatrix{
H^1(Y',\Z/p\Z) \ar@{^{(}->}[r] \ar[d]^{\Psi_{Y,y,\Z/p\Z}} & 
H^1(Y',G)
\ar[d]^{\Psi_{Y,y,G}} \\
H^1(k((t)),\Z/p\Z) \ar@{^{(}->}[r]  & 
H^1(k((t)),G)
}
\]

and since the top and right arrows are injective, we know that $\Psi_{Y,y,\Z/p\Z}$ is injective.

Next we show that (3) implies (1), so assume that $\Psi_{Y,y,\Z/p\Z}$ is injective. Let $b$ be an element of $k((t))$ such that $\wp(b) \in \mathcal{O}(Y')$. 
Then, by Lemma~\ref{mult}, we have that $k((t))[x]/(x^p - x)$ is isomorphic to $k((t))[x]/(x^p - x - \wp(b))$ as $\Z/p\Z$-Galois $k((t))$-algebras. 
But since $\Psi_{Y,y,\Z/p\Z}$ is injective, we have that $\mathcal{O}(Y')[x]/(x^p -x) \cong \mathcal{O}(Y')[x]/(x^p -x-\wp(b))$, so by Lemma~\ref{mult}, $0 = \wp(c)+\wp(b)$ for some $c \in \mathcal{O}(Y')$, which means $b - c \in \mathbb{F}_p$, so $b \in \mathcal{O}(Y')$, and so $\wp(k((t))) \cap \mathcal{O}(Y') \subseteq \wp(\mathcal{O}(Y'))$. The direction $\wp(k((t))) \cap \mathcal{O}(Y') \supseteq \wp(\mathcal{O}(Y'))$ is clear. \qed

\subsection{Proof of Theorem \ref{surj}}
We first note that (1) implies (3) by Artin--Schreier (a special case of Lemma~\ref{mult}).

Next, we show that (3) implies (2), so suppose $\Psi_{Y,y,\Z/p\Z}$ is surjective for all $n$. Since $U_n(\mathbb{F}_p)\cong \Z/p\Z$,
we can also assume that $\Psi_{Y,y,U_n(\Z/p\Z)}$ is surjective. We proceed by induction on the order of $G$. Since $G$ is a $p$-group, $G$ has a central subgroup $H$ which is isomorphic to $\Z/p\Z$.

By Lemma~\ref{longexact} of this paper, and by
Lemma~1.4.3 of \cite{katz}
(the latter stating that both $H^2(k((t)), \Z/p\Z)$ and $H^2(Y', \Z/p\Z)$ are zero), 
the map $\Spec k((t)) \to Y'$ induces the following commutative diagram
\[
\xymatrix{
1 \ar[r]& 
H^1(Y',\Z/p\Z)\ar[r]^{\iota_{Y'}} \ar[d]^{\Psi_{Y,y,\Z/p\Z}} & 
H^1(Y',G) \ar[r]^{\phi_{Y'}} \ar[d]^{\Psi_{Y,y,G}} & 
H^1(Y',G/H) \ar[r] \ar[d]^{\Psi_{Y,y,G/H}} & 
1\\
1 \ar[r] & 
H^1(k((t)),\Z/p\Z) \ar[r]^{\iota_{k((t))}}  & 
H^1(k((t)),G) \ar[r]^{\phi_{k((t))}}  & 
H^1(k((t)),G/H) \ar[r]  & 
1
}
\]
whose rows are exact. 
Also by Lemma~\ref{longexact}, 
two elements of $H^1(Y', G)$ have the same image in $H^1(Y', G/H)$ if and only if they are in the same $H^1(Y', \Z/p\Z)$-orbit (and similarly for $H^1(k((t)), G))$.
And by the assumption (3) and the inductive hypothesis, $\Psi_{Y,y,\Z/p\Z}$ and $\Psi_{Y,y,G/H}$ are surjective. The surjectivity of $\Psi_{Y,y,G}$ is proved via the following diagram chase.

Let
$\tilde{\beta}$ be an element of $H^1(k((t)),G)$, and let $\tilde{\gamma}:= \phi_{k((t))}(\tilde{\beta})$. By the inductive hypothesis, there exists $\gamma \in H^1(Y',G/H)$ such that $\Psi_{Y,y,G/H}(\gamma) = \tilde{\gamma}$.
Let $\beta$ be an element of $H^1(Y',G)$ mapping to $\gamma$. 
Then there exists
$\tilde{\alpha} \in H^1(k((t)), \Z/p\Z)$ such that $\iota_{k((t))} (\tilde{\alpha} ) \cdot \Psi_{Y,y,G}(\beta) = \tilde{\beta}$, and by the inductive hypothesis there is an element $\alpha \in H^1(Y', \Z/p\Z)$ mapping to $\tilde{\alpha}.$ Then $\Psi_{Y,y,G}(\alpha \cdot \beta) = \tilde{\beta}$, so $\Psi_{Y,y,G}$ is surjective.

Next, we show that (2) implies (1), so suppose $\Psi_{Y,y,G}$ is surjective for some finite $p$-group $G$, and again let $H$ be a nontrivial central subgroup of $G$ isomorphic to $\Z/p\Z$. The diagram above shows that $\Psi_{Y,y,G/H}$ is surjective; iterating this process shows that $\Psi_{Y,y,\Z/p\Z}$ is surjective.
Let $f \in k((t))$, so $k((t))[x]/(x^p-x-f)$ is a $\Z/p\Z$-Galois \'{e}tale algebra over $k((t))$.
Since $\Psi_{Y,y,\Z/p\Z}$ is surjective, Theorem~\ref{inaba} tells us that $f$ is $p$-equivalent to an element $g$ of $\mathcal{O}(Y')$, so $f = \wp(b) + g$ for some $b \in k((t))$. Thus the equality in (1) holds. \qed

When $\Psi_{Y,y,G}$ is an injection, we now use the Inaba classification to give a concrete description of the unique cover mapping under $\Psi_{Y,y,G}$ to a given $G$-Galois cover of $k((t))$ which is in the image of $\Psi_{Y,y,G}$. We use the notation and index-ordering described in Theorem~\ref{inaba}.

%We now relate the structure of the subextentions described in part (3) of Theorem~\ref{inaba} over $\mathcal{O}(Y')$ and $k((t))$ when $\psi_{Y,y,G}$ is an injection, using the notation and index-ordering described above Theorem~\ref{inaba}, and connect this to equivalence of categories. In the following lemma, we use the Inaba classification to give a concrete construction of the unique cover mapping to a $G$-Galois field extension.

\begin{lem}\label{subextensions}
Suppose that $\Psi_{Y, y, G}$ is an injection, and that  $L$ is a $G$-Galois \'{e}tale $k((t))$-algebra
%; by Theorem~\ref{inaba} and Theorem~\ref{surj}, %
which is isomorphic to
%$L \cong 
$k((t))[A]$ with $A^{(p)} = MA$ and $M \in  U_n(\mathcal{O}(Y'))$. 
Let $E$ denote the subring of $L$ generated by $\mathcal{O}(Y')$ and the entries $a_{ij}$ of $A$, and let $E_{i,j} = \mathcal{O}(Y') [a_{kl}: (k,l)<(i,j)]$. 
Then $L_{i,j}/L_{<i,j}$ is a nontrivial extension if and only if $E_{i,j}/E_{<i,j}$ is a nontrivial extension, and $\Psi_{Y,y,G}(E) = L$.
\end{lem}
\begin{proof}
Since $M\in U_n(\mathcal{O}(Y'))$, we have that $E_{i,j} = E_{<i,j}[a_{ij}]$ where $\wp(a_{ij}) \in E_{<ij}$.

If $E_{i,j}/E_{<i,j}$ is a trivial extension, then $a_{ij} \in E_{<ij} \subset L_{<ij}$, so $L_{i,j}/L_{<i,j}$ is a trivial extension. 

We next show that if $L_{i,j}/L_{<i,j}$ is a trivial extension, then $E_{i,j}/E_{<i,j}$ is trivial. So suppose that $a_{ij} \in L_{<i,j}$. We know that $\wp(a_{ij}) \in E_{<i,j}$ and we want to show that $a_{ij}$ is also in $E_{<i,j}$. More generally we show by induction that for any $(k,\ell)$ and $a\in L_{k,l}$, if $\wp(a) \in E_{k,l}$, then $a\in E_{k,l}$. 
The base case is that $\wp(k((t))) \cap \mathcal{O}(Y') = \wp(\mathcal{O}(Y'))$; this follows from Theorem~\ref{inj}.
So suppose $L_{k,l}$ strictly contains $k((t))$.
We may assume  $L_{k,l}$ is isomorphic to  $L_{<k,l}[x]/(x^p-x-f)$ for an indeterminate $x$ and $f \in L_{<k,l}$, and so $a$ can be expressed uniquely as a sum
\[
a = \sum_{m=0}^{p-1} c_m x^m
\]
where $c_m \in L_{<k,j}$.
The $p-1$st coefficient of $\wp(a)$ is then $\wp(c_{p-1})$, which is in $E_{<k,l}$, and hence by induction $c_{p-1}$ is in $E_{<k,l}$.
Subtracting off the top-degree term of $a$ and iterating this procedure, we can conclude that all coefficients of $a$ are in $E_{<k,l}$,
and since $x \in E_{k,l}$, we conclude that $a \in E_{k,l}$, which is what we wanted to show.

Since $E$ is described as an iterative extension in a completely analogous way to $L$ as in Theorem~\ref{inaba}, it follows that $E$ is an $\et$ $G$-Galois $\mathcal{O}(Y')$-algebra with $L = \Psi_{Y,y,G}(E)$.
\end{proof}

\begin{thm}\label{equivcat}
For a nontrivial $p$-group $G$, the map $\Psi_{Y,y,G}$ is a bijection if and only if $\psi_{Y,y,G}$ is an equivalence of categories.
\end{thm}
\begin{proof}
First suppose that $\psi_{Y,y,G}$ is an equivalence of categories. Then since  $\psi_{Y,y,G}$ is essentially surjective,  $\Psi_{Y,y,G}$ is surjective. And if $\phi: L_1 \to L_2$ is an isomorphism of $G$-Galois $k((t))$-algebras and $\Psi_{Y,y,G}(E_i) = L_i$ for $i=1,2$, then since $\psi_{Y,y,G}$ is fully faithful, $\phi$ lifts to an isomorphism $E_1 \to E_2$, so $\Psi_{Y,y,G}$ is injective.

For the converse, suppose that $\Psi_{Y,y,G}$ is a bijection. 
%This implies that we have nicely corresponding decompositions of $k((t))$- and $\mathcal{O}(Y')$-algebras as detailed in Lemma~\ref{subextensions}.
We construct an inverse functor to $\psi_{Y,y,G}$.
Let $L$ be a $G$-Galois \'{e}tale $k((t))$-algebra, so $L$ is of the form $k((t))[A]$ where $A^{(p)} = MA$ for some $M\in U_n(\mathcal{O}(Y'))$ as argued in the proof of Lemma~\ref{subextensions}.
We define $\psi_{Y,y,G}^{-1}(L)$ to be $E$ as defined in Lemma~\ref{subextensions}.
Now let $L=k((t))[A]$ and $L'=k((t))[A']$ be two $G$-Galois \'{e}tale $k((t))$-algebras, with $A^{(p)} = MA$ and $A'^{(p)} = M'A'$, for $M,M'\in U_n(\mathcal{O}(Y'))$. Suppose $\phi$ is a morphism $L\to L'$; we will show that $\phi$ restricts to a homomorphism $E\to E'$.
We induct on $(i,j)$, showing that $\phi(E_{i,j})\subset E'$ for all $(i,j)$.
The base case follows from $\phi$ being a map of $k((t))$-algebras and hence a map of $\mathcal{O}(Y')$-algebras.
So suppose that $\phi(E_{<i,j})\subset E'$.
We have that $E_{i,j}=E_{<i,j}[a_{ij}]$, so we just need to show that $\phi(a_{ij}) \in E'$. We know that $\wp(a_{ij})$ is in $E_{<i,j}$, and so $\wp(\phi(a_{ij}))=\phi(\wp(a_{ij}))$ is in $E_{<i,j}$ by induction. But then $\phi(a_{ij}) \in E_{<i,j}$ by an argument as in the proof of Lemma~\ref{subextensions}. So we define $\psi_{Y,y,G}^{-1}(\phi)$ to be this restriction.
Then it is straightforward to see that $\psi_{Y,y,G}^{-1}$ and $\psi_{Y,y,G}$ are inverses.
\end{proof}

We now give a concise reformulation of the criterion for when $\Psi_{Y,y,G}$ is a bijection.
We denote by $F$ the Frobenius map $\mathcal{O}_Y \to \mathcal{O}_Y $ sending $f \mapsto f^p$. This induces a map $F^* : H^1(Y,\mathcal{O}_Y) \to H^1(Y,\mathcal{O}_Y )$; we also let $\wp^*:= F^* - Id$ be the map on cohomology induced by Artin--Schreier. Then our reformulation is as follows.

\begin{cor}\label{bij}
The map $\psi_{Y,y,G}$ is an equivalence of categories if and only if $\wp^*: H^1(Y,\mathcal{O}_Y) \to H^1(Y,\mathcal{O}_Y)$ is a bijection.
\end{cor}

\begin{proof}
The natural open immersion $i: Y' \hookrightarrow Y$ gives the following exact sequence of sheaves on $Y$:
\[
0 \to 
\mathcal{O}_Y \to 
i_* \mathcal{O}_{Y'} \to 
\textrm{skysc}_y \left( \frac{k((t))}{k[[t]]} \right) \to
0
\]
where $\textrm{skysc}_y \left( \frac{k((t))}{k[[t]]} \right)$ denotes the skyscraper sheaf at $y$, with value group $k((t))/k[[t]]$ where the group structure is given by the additive structure on $k((t))$.
We see that $H^1(Y,i_*\mathcal{O}_{Y'}) \cong H^1(Y',\mathcal{O}_{Y'})= 0$ 
since $Y'$ is affine, 
 so the induced long exact sequence in Zariski cohomology gives us an isomorphism
\[
\frac{
H^0
\left(
Y,
\textrm{skysc}_y \left( \frac{k((t))}{k[[t]]} \right)
\right)
}{
H^0(Y,i_*\mathcal{O}_{Y'})
}
=
\frac{
k((t))
}{
\mathcal{O}(Y') + k[[t]]
}
\xrightarrow{\sim}
H^1(Y,\mathcal{O}_Y).
\]
We also see that the map
\[
F^* :
\frac{
k((t))
}{
\mathcal{O}(Y') + k[[t]]
}
\to 
\frac{
k((t))
}{
\mathcal{O}(Y') + k[[t]]
}
\]
maps $\overline{f}\mapsto \overline{f}^p$ for $f \in k((t))$.

Now suppose $\wp^*$ is surjective. Then for any $f \in k((t))$ there exist $g \in k((t))$, $h \in \mathcal{O}(Y')$, $l \in k[[t]]$ such that $f = g^p - g + h + l$.
Let $a_0$ be the constant term of $l$; by setting $h' := h + a_0$ and $l' := l - a_0$ we can assume $l$ has no constant term. 
So $l = \sum_{i=1}^\infty a_i t^i$.
We define a power series $\tilde{l}:= \sum_{i=1}^\infty b_i t^i$ where $b_i = -a_i$ for $i$ not divisible by $p$ and $b_{np} = b_n^p - a_{np}$.
So $f = \wp(g + \tilde{l}) + h$ and $\Psi_{Y,y,G}$ is surjective for all nontrivial $p$-groups $G$ by Theorem~\ref{surj}. The converse, that surjectivity of the maps $\Psi_{Y, y, G}$ implies surjectivity of $\wp^*$, is straightforward.

Next, suppose $\wp^*$ is injective, and consider $f$ such that $\wp(f) \in \mathcal{O}(Y')$.
Then $f \in k[[t]] + \mathcal{O}(Y')$,
so there exist $g \in \mathcal{O}(Y')$, $l \in k[[t]]$ such that $f = g + l$. Again, we can assume that $l$ has no constant term, so $l \in tk[[t]]$.
But then $\wp(g)  + \wp (l) \in \mathcal{O}(Y')$, which implies $\wp(l) \in \mathcal{O}(Y')$. 
But since a nonzero $l \in tk[[t]]$ would have a zero at $y$, it could not come from a regular function on $Y'$, so we must have $l=0$ and so in fact $\wp (f) \in \wp(\mathcal{O}(Y'))$ and $\Psi_{Y,y,G}$ is injective for all nontrivial $p$-groups $G$ by Theorem~\ref{inj}.
Conversely, suppose $\Psi_{Y,y,G}$ is injective for all nontrivial $p$-groups $G$, and consider $f \in k((t))$ such that $\wp(f) = g + l$ for some $g \in \mathcal{O}(Y')$, $l \in tk[[t]]$.
We can write $l = \wp (\tilde{l})$ for $\tilde{l}$ as above, so $\wp(f-\tilde{l}) \in \mathcal{O}(Y')$, which implies $f-\tilde{l} \in \mathcal{O}(Y')$ by the hypothesis.
Then $\overline{f} = \overline{0}$ in $k((t))/(k[[t]] + \mathcal{O}(Y'))$, which is what we wanted to show.
\end{proof}

\section{Examples}
In this section, we apply the previous results to study $\psi_{E,O}$ for $E$ an elliptic curve over $\mathbb{F}_p$ and $O$ the point at infinity. The operator $F^*$ on $H^1(E,\mathcal{O}_E)$ acts as multiplication by some $a \in \mathbb{F}_p$, and in fact $\# E(\mathbb{F}_p) = p+1-a$. We say that $E$ is anomalous if $\# E(\mathbb{F}_p) = p$.

\begin{thm}\label{elliptic}
For an elliptic curve $E$ over $\mathbb{F}_p$, the following are equivalent:
\begin{enumerate}
    \item $E$ is not anomalous.
    \item The map $\Psi_{E,O,G}$ is injective for every nontrivial $p$-group $G$.
    \item The map $\Psi_{E,O,G}$ is surjective for every nontrivial $p$-group $G$.
    \item The map $\psi_{E,O}$ gives an equivalence of categories.
\end{enumerate}
\end{thm}

\begin{proof}
We let $e$ denote a generator of $H^1(E,\mathcal{O}_E)$ as an $\mathbb{F}_p$-vector space, 
so $F^*$ maps $e \mapsto a\cdot e$. 
Then for $x\in \mathbb{F}_p$, we have $xe \mapsto x^p a e$. 
But since $x\in \mathbb{F}_p$, the map is $xe\mapsto xae$,
so $\wp^*(xe)=x(a-1)e$. Then the map 
$\wp^* : H^1(E,\mathcal{O}_E) \to H^1(E,\mathcal{O}_E)$ is surjective if and only if $a\neq 1$
 and it’s injective if and only if $a\neq  1$. Applying Corollary~\ref{bij}
 gives the result.
\end{proof}

As indicated in \cite{olson}
anomalous curves are, as their name suggests, uncommon, and so Theorem~\ref{elliptic}
provides us with a broad class of curves whose $p$-group Galois covers correspond nicely to $p$-group $k((t))$-extensions.

\bibliographystyle{amsalpha}
\bibliography{LTG}

\end{document}